\documentclass[11pt]{amsart}

\usepackage{amsfonts, enumerate} 
\usepackage{amssymb}
\usepackage[final]{graphicx, graphics,stmaryrd} 
\usepackage[margin=3cm]{geometry} 
\usepackage{amscd}
\usepackage[all]{xy}
\usepackage[breaklinks,bookmarksopen,bookmarksnumbered]{hyperref}

\usepackage{tensor}
\usepackage[english]{babel}
\usepackage[utf8]{inputenc}
\usepackage{amsthm}
\usepackage{graphics}
\usepackage{amsmath}
\usepackage{amstext}
\usepackage{engrec}
\usepackage{rotating}
\usepackage[safe,extra]{tipa}
\usepackage{multirow}
\usepackage{enumitem}
\usepackage{bm}
\usepackage{stackrel} 
\usepackage{tikz-cd} 
\usepackage{mathtools} 
\numberwithin{equation}{section} 

\newtheorem{teo}{Theorem}[section]

\newtheorem{prop}[teo]{Proposition}
\newtheorem{lemma}[teo]{Lemma}

\theoremstyle{definition}
\newtheorem{defin}[teo]{Definition}
\newtheorem{rmk}[teo]{Remark}

\newcommand{\mc}{\mathcal}

\def\R{\mathbb{R}}

\def\N{\mathbb N}

\def\d{\mathrm{d}}

\newcommand{\m}{\mbox}

\newcommand{\cor}{\textit}

\newcommand{\fine}{\qed\newline}
\newcommand{\xx}{\otimes}

\DeclareMathOperator{\id}{id}

\DeclareMathOperator{\Div}{Div}

\DeclareMathOperator{\supp}{supp}

\def\Y{\mc{Y}}
\newcommand{\de}{\partial}

\DeclareMathOperator{\Geod}{Geod}

\DeclareMathOperator{\Kan}{Kan}
\DeclareMathOperator{\Mon}{Mon}
\DeclareMathOperator{\BB}{BB}

\newcommand\norm[1]{\lVert#1\rVert}

\begin{document}
		\title[Benamou--Brenier and Kantorovich on sub-Riemannian manifolds]{Benamou--Brenier and Kantorovich on sub-Riemannian manifolds with no abnormal geodesics}
	
	\author[G.\,Citti]{Giovanna Citti}
	\author[M.\,Galeotti]{Mattia Galeotti}
	\author[A.\,Pinamonti]{Andrea Pinamonti}
	
		\address[G.\,Citti]{Dipartimento di Matematica
		\newline\indent Università  di Bologna \newline\indent
		Piazza di Porta San Donato 5, 40122 Bologna, Italy}
	\email{giovanna.citti@unibo.it}
	
	\address[M.\,Galeotti]{Dipartimento di Matematica
		\newline\indent Università  di Bologna \newline\indent
		Piazza di Porta San Donato 5, 40122 Bologna, Italy}
	\email{mattia.galeotti4@unibo.it}
	
		\address[A.\,Pinamonti]{Dipartimento di Matematica
		\newline\indent Università degli Studi di Trento \newline\indent
		Via Sommarive, 14, 38123 Povo, Trento, Italy}
	\email{andrea.pinamonti@unitn.it}
	
%

	

\maketitle
\begin{abstract}
We prove that the Benamou--Brenier formulation of the Optimal Transport
problem and the Kantorovich formulation are equivalent on a sub-Riemannian connected and
complete manifold $M$ without boundary and with no non-trivial abnormal geodesics, when the problems
are considered
between two measures with finite $2$-momentum.
Furthermore, we prove the existence of a minimizer for the Benamou--Brenier
formulation and link it to the optimal transport plan.
\end{abstract}

\section{Introduction}

{
The Benamou--Brenier ``dynamical'' formulation of Optimal Transport
consists in the minimization problem for the functional
\begin{equation}\label{bbfunctional}
	J_{\BB}(\mu_t,v_t)=\int_0^1 \|v_t(x)\|^2 \, \mathrm{d}\mu_t(x)\, \mathrm{d}t,
\end{equation}
where \(v_t\) is a Borel family of vector fields and \(\mu_t\) is a narrowly continuous
family of probability measures connecting the prescribed marginals \(\mu_0\) and \(\mu_1\).
Both \(\mu_t\) and \(v_t\) are defined on a differentiable manifold \(M\)
and are constrained by the \emph{Continuity Equation}
\begin{equation}\label{cequation}
	\dot\mu_t+\Div(\mu_t v_t)=0.
\end{equation}
For a detailed discussion of the Continuity Equation we refer to
\cite[Chapter~8]{ags08} and \cite{amcri14}.
This equation links the curve \(\mu_t\) to the flow of the ODE associated with \(v_t\),
encoding the conservation of mass for the action~\eqref{bbfunctional}.
Altogether, this provides a time-dependent approach on the transport problem. 

By contrast, in the classical Kantorovich formulation
the dynamical component is encoded in the squared distance function,
\[
 d^2(x,y):=\inf\left\{\left.\int_0^1 \|\dot\omega_t\|^2\,\mathrm{d}t\ \right|\ 
 \omega\colon [0,1]\to M\ \text{absolutely continuous},\,
 \omega_0=x,\,\omega_1=y\right\}.
\]
This is the so-called \emph{cost function}, and the corresponding minimization problem concerns
the coupling of points. In this ``static'' formulation one seeks to minimize
\begin{equation}\label{kanfunctional}
	J_{\Kan}(\gamma)=\int_{M\times M} d^2(x,y)\,\mathrm{d}\gamma(x,y),
\end{equation}
where \(\gamma\) ranges among admissible transport plans, that is,
probability measures on \(M\times M\) whose marginals coincide with
\(\mu_0\) and \(\mu_1\) on the two coordinate projections.

The Kantorovich and Benamou--Brenier formulations
have been shown to be equivalent on bounded Euclidean domains \cite{benabre00},
and more generally on complete Riemannian manifolds without boundary,
see, for instance, \cite{amgi13}.
In the present work we establish their equivalence on a sub-Riemannian manifold \(M\),
under a suitable set of assumptions specified below.\newline

The transport problem was first introduced by Monge at the end of the 18th century.
It consists in minimizing the functional
\[
	J_{\Mon}(T)=\int_M c(x,T(x))\,\mathrm{d}x
\]
for a cost function \(c\colon M\times M\to \mathbb{R}_{\geq 0}\cup\{+\infty\}\),
where \(T\colon M\to M\) satisfies the transport constraint \(T_\#\mu_0=\mu_1\).
However, this formulation is ill-suited for general analysis,
as the transport constraint is not closed with respect to the natural weak topologies.

The Kantorovich relaxation overcomes this difficulty.
In the Kantorovich setting, the existence of a minimizer (an optimal plan)
is guaranteed for any Polish space \((M,d)\)
whenever the cost function is lower semicontinuous and bounded from below,
see \cite[Theorem~1.5]{amgi13} and also \cite[Theorem~5.10]{villa08} for a broader discussion.
Moreover, in the Riemannian case with cost given by the squared distance,
one also obtains the existence and uniqueness of an optimal map,
that is, a measurable map inducing an optimal plan.
This was established by Brenier in the Euclidean case \cite{bre91}
and by McCann in the general Riemannian setting \cite{mccann01}.\newline

In the first decades of the 21st century, renewed interest arose in Optimal Transport
within the framework of sub-Riemannian geometry.
A sub-Riemannian manifold is a differentiable manifold endowed with
a totally non-integrable distribution \(HM\subset TM\)
and a positive definite quadratic form \(g_{SR}\) defined on \(HM\),
see \cite{riffo14} for an introduction.
Several existence, uniqueness, and regularity results for optimal maps
have been obtained in particular cases, including the Heisenberg group \cite{amrigo04},
\(2\)-generating distributions \cite{agralee09}, and more general settings
\cite{figajui08, figariffo10}. However, the general problem remains open.

While in the Riemannian case all geodesics are \emph{normal},
in the sub-Riemannian framework a second type, the so-called \emph{abnormal} geodesics, may occur.
The presence of such geodesics causes the sub-Riemannian distance
to lose local convexity, which constitutes the main obstacle
to the existence of optimal maps.\newline

To address the dynamical problem in the sub-Riemannian setting,
we consider a \emph{relaxed version} of the Benamou--Brenier problem.
Specifically, we minimize the integral of \(\|v\|^2\) with respect to
a probability measure \(\eta\) on \([0,1]\times HM\),
instead of integrating a single vector field \(v_t\) against \(\mu_t\) at each time.
This approach, already employed in \cite{ello23, elam24} for \(M=\mathbb{R}^n\)
with a sub-Riemannian structure, fits naturally within the framework of
control relaxation techniques (see also
\cite{buttadalma82, butta87, masmig89, prr20, savasodi24}).

For a rigorous formulation of the relaxed problem we employ the class of Young measures.
Given a Polish space \((X,d)\), a Young measure \(\eta\)
is a positive measure on \([0,1]\times X\)
whose first-coordinate projection \(\pi^{(t)}_\#\eta\) is Lebesgue  measurable,
hence a disintegration \(\eta_t\) exists with respect to Lebesgue measure on \([0,1]\).
If \(\pi\colon TM\to M\) denotes the canonical projection,
then the pushforward of a Young probability on \([0,1]\times HM\)
yields a narrowly continuous curve of probability measures
\(\mu_t=\pi_\#\eta_t\in\mathcal{P}(M)\).

The Continuity Equation can be adapted to this framework (see Equation~\eqref{transport})
and used as a constraint for the minimization of the relaxed functional
\begin{equation}\label{relaxedfunctional}
	J_{\BB}^\star(\eta)=\int_{[0,1]\times HM}\|v\|^2\,\mathrm{d}\eta(t,v),
\end{equation}
where we also require that \(\pi_\#\eta_0=\mu_0\) and \(\pi_\#\eta_1=\mu_1\).
A Young measure $\eta$ on $I\times HM$ respecting
the Continuity Equation is called a \cor{transprot measure}.\newline

We first prove the equivalence between the original (sub-Riemannian)
Benamou--Brenier formulation and its relaxed version.
Then we show the equivalence with the Kantorovich formulation
by means of a \emph{superposition principle}.
As explained in detail in \S\ref{superyoung},
for any transport measure $\eta$ there exists a probability
 $\widetilde \eta$ on the space of (generalized) curves on $M$,
such that if $\overline e_t$ is the map sending a (generalized) curve $\varphi$
to the point $\varphi(t)$, then
\[
\eta=(\overline e_t)_\#\widetilde \eta.
\]

\begin{teo}\label{mainteo}
	Let \(M\) be a connected, complete sub-Riemannian manifold
	without boundary and without non-trivial abnormal geodesics.
	Let \(\mu_0,\mu_1\in\mathcal{P}_2(M)\).
	Then the infima of the Kantorovich, Benamou--Brenier,
	and relaxed Benamou--Brenier transport problems are all finite,
	attained, and equal. If $\gamma$ is a minimizer
	for the Kantorovich problem,
	then the pushforward $F_\#(\mc L\xx \gamma)$ is a minimizer
	for the Benamou--Brenier problem, where $\mc L$ is the Lebesuge measure
	and $F$ the map defined in \eqref{Fmap}.
	
	Finally, for any Benamou--Brenier minizimer $\eta$,
	its decomposition $\widetilde \eta$ is supported
	on (generalized and constant speed) geodesics.
\end{teo}

Compared with \cite{elam24}, our approach is developed
outside the control-theoretic framework.
In addition, with respect to the Euclidean setting,
a further difficulties arises. 
The use of the superposition principle depends
on the construction of a measurable map \(S\)
that associates to each pair \((x,y)\in M\times M\)
a geodesic joining \(x\) to \(y\).
The construction of such a map involves nontrivial technical challenges,
which we overcome using Suslin set theory, as shown in Theorem~\ref{teoSuslin}.

The main obstacle to extending our results to manifolds
admitting non-trivial abnormal geodesics lies precisely
in the measurable construction of the map \(S\).\newline

The paper is organized as follows.
Sections~\S\ref{srstruc} and~\S\ref{otinsr} recall
the basic notions on sub-Riemannian manifolds and Optimal Transport.
Section~\S\ref{relax} is devoted to the formulation of the relaxed problem
and to the proof of its equivalence with the standard Benamou--Brenier formulation,
while Section~\S\ref{kanequiv} focuses on the equivalence between
the Kantorovich and Benamou--Brenier formulations.
The Appendix collects several technical results used throughout the paper.
In particular, \S\ref{transpmeas} recalls the main results on Young transport measures,
\S\ref{Smap} establishes the existence of the measurable map
\(S\colon M\times M\to \mathrm{Geod}(M)\),
\S\ref{superyoung} presents and proves the superposition principle
in the sub-Riemannian setting, and
\S\ref{cptresult} presents a compactness result
for the sublevel sets of the relaxed functional.\newline
}

\section{Preliminaries}

In this section we introduce the basics
of sub-Riemannian structures on manifolds,
giving in particular some general results on the sub-Riemannian
distance and sub-Riemannian geodesics.
Furthermore,
we introduce optimal transport on metric
spaces and more specifically the formulations of Kantorovich
and Benamou--Brenier in the context of sub-Riemannian manifolds.

\subsection{Sub-Riemannian structures}\label{srstruc}

A sub-Riemannian manifold
is a triplet $(M,H,g_{SR})$ of a differentiable manifold $M$,
a completely non-integrable distribution $HM\subset TM$
and a positive smooth quadratic form $g_{SR}$ defined on $HM$. By completely
intregable we mean that for every $x\in M$, if $X_1,\dots,X_m$ is a base of $HM$
on an open neighborhood $U\subset M$ such that $x\in U$, then there exists $r\in \N^+$
such that the brackets 
\[
\left[X_{i_1},[\dots,[X_{i_{k-1}},X_{i_k}]]\right](x)\quad \m{with }i_1,\dots,i_k\in\{1,\dots,m\}\m{ and }1\leq k\leq r,
\]
span the whole tangent space $T_xM$.	

For a wide introduction to sub-Riemannian manifolds and their metric and differential structure,
see \cite{abb19} and \cite{bella96}. Here, we recall the basics of the metric
structure induced by the sub-Riemannian setting.\newline

We denote by $I$ the unitary interval $[0,1]$.
We call \cor{admissible} any absolutely continuous curve $\omega\colon I\to M$
such that the tangent vector $\dot \omega_t$ lies in $H_{\omega_t}M$ for a.e.~$t$
in $I$. These are also called \cor{horizontal} curves.
We denote by 
\[
\Omega_{HM}:=\left\{\omega\colon I\to M:\ \m{absolutely continuous and }\dot\omega_t\in H_{\omega_t}M\ \m{for a.e.}\,t\in I\right\}
\]
the space of admissible curves.
A sub-Riemannian distance $d_{SR}$ is defined on $M$
as 
the infimum of the lenght of admissible curves
between two points,
\[
d_{SR}(x,y):=\inf \left\{\int_I||\dot\omega_t||\,\d t:\ \omega\in \Omega_{HM},\ \omega_0=x,\ \omega_1=y\right\}.
\]
The celebrated Chow–Rashevskii theorem states that if $M$ is connected,
then $d_{SR}$ is always finite, \cor{i.e.}~for any pair $(x,y)$ there exists an
admissible curve joining $x$ to $y$.
Moreover, if $(M,d_{SR})$ is complete, for any pair $(x,y)$ there exists
a length minimizer joining $x$ to $y$, \cor{i.e.}~the distance $d_{SR}(x,y)$ is effectively realized.

When talking about geodesics we consider constant speed geodesics, 
the admissible curves with constant speed that minimize the length functional.
Any geodesic can be reparametrized to a constant speed geodesic.

Equivalently, constant speed geodesics minimize the energy functional,
\[
\int_I ||\dot\omega_t||^2\,\d t,
\]
when the extremes of the curve $\omega_t$ are fixed. Therefore,
\[
\Geod(M):=\left\{\omega\colon I\to M:\ \m{absolutely continuous and }d^2_{SR}(\omega_0,\omega_1)=\int_I||\dot\omega_t||^2\d t\right\}.
\]
We consider the space of (constant speed) geodesics as endowed with the sup distance
\[
d(\omega^{(1)},\omega^{(2)}):=\sup_{t\in I}d_{SR}(\omega^{(1)}(t),\omega^{(2)}(t))\quad\forall \omega^{(1)},\omega^{(2)}\in \Geod(M).
\]\newline

Any geodesic $\omega_t$ on a sub-Riemannian manifold corresponds to the projection
of an integral curve $\lambda_t\colon I \to T^*M$ of the Hamiltonian vector field
$X_H\in TT^*M$. Any such curve $\lambda$ satisfies the Hamilton equation
and the value of the Hamiltonian $H(\lambda_t)$ is constant for $t\in I$.
From this characterization two types of geodesics arise:
\begin{enumerate}
	\item \cor{normal} geodesics are the projection of curves with $H(\lambda_t)>0$;
	\item \cor{abnormal} geodesics are the projection of curves with $H(\lambda_t)=0$.
\end{enumerate}
A geodesic is called {strictly normal} (resp.~{strictly abnormal}) when it is normal and not abnormal (resp.~abnormal and not normal).
Normal geodesics are uniquely determined by the initial conditions of the associated curve $\lambda_t\in T^*M$,
namely the initial covector $\lambda_0\in T^*_xM$.\newline

In the following we work with a complete connected sub-Riemannian manifold $M$
without boundary and with no non-trivial abnormal geodesics. As there is no risk
of confusion, from now on we  use the simpler notation $d(\cdot,\cdot)$ for
the sub-Riemannian distance $d_{SR}(\cdot,\cdot)$.\newline

\subsection{Optimal transport on sub-Riemannian manifolds}\label{otinsr}

The Optimal Transport problem may be introduced in full generality
over a Polish space $(X,d)$.
For a wide overview of this subject see for example \cite{villa08, amgi13}.
We use the notation $\mc P(X)$ for the space of Borel probabilty measures on $X$
and, given a Borel map $f\colon X\to Y$ between Polish spaces, we denote by
$T_\#\mu$ the pushforward of $\mu\in \mc P(X)$, 
\[
T_\#\mu(B)=\mu(T^{-1}(B)),\ \ \forall B\in \mc B(Y).
\]

Over the probability space $\mc P(X)$ we will often consider
convergence with respect to the narrow topology,
meaning that $\mu_n\to \mu$ if 
\[
\int_X\phi\,\d\mu_n\to\int_X\phi\,\d\mu\ \ \ \forall \phi\in C_b(X).
\]\newline

Once we fix a (Borel) cost function $c\colon X\times Y\to \R\cup\{+\infty\}$,
we can introduce the Monge version of the transport problem.
Given $\mu_0\in \mc P(X),\ \mu_1\in \mc P(Y)$, we want to minimize
\begin{equation}\tag{Mon}\label{monge}
	C_{\Mon}(\mu_0,\mu_1):=\inf_{T_\#\mu_0=\mu_1}\int_X c(x,T(x))\,\d\mu_0(x)
\end{equation}
where the infimum is taken among all Borel maps $T\colon X\to Y$ such that $T_\#\mu_0=\mu_1$.
This formulation has a series of patologies, for instance
the constraint $T_\#\mu_0=\mu_1$ lack closedness with respect to the principal
weak topologies.

The problem was formulated in a more general fashion by Kantorovich.
We denote by $\Pi(\mu_0,\mu_1)$ the set of \cor{admissible transport plans} between $\mu_0$ and $\mu_1$,
\[
\Pi(\mu_0,\mu_1):=\left\{\gamma\in \mc P(X\times Y):\ \pi^{(X)}_\#\gamma=\mu_0,\ \pi^{(Y)}_\#\gamma=\mu_1\right\}
\]
where $\pi^{(X)}$ and $\pi^{(Y)}$ are the projections from $X\times Y$ to $X$ and $Y$ respectively.
The minimization problem becomes
\begin{equation}\tag{Kan}\label{kantorovich}
	\begin{split}
		J_{\Kan}(\gamma)&:=\int_{X\times Y}c(x,y)\d\gamma(x,y)\\
		\\
	C_{\Kan}(\mu_0,\mu_1)&:=\inf_{\gamma\in \Pi(\mu_0,\mu_1)}J_{\Kan}(\gamma).
	\end{split}
\end{equation}
Observe that a transport map $T$,
as the one considered in the Monge formulation,
is represented in the Kantorovich formulation
via the induced plan $\gamma=(\id_X\times T)_\#\mu_0$,
therefore \eqref{kantorovich} is a generalization of~\eqref{monge}.
Moreover, $\Pi(\mu_0,\mu_1)$ is always non-empty, convex and compact with respect
to the narrow topology. A minimizer for \eqref{kantorovich} always exists under
very general assumptions.
\begin{teo}[See Theorem 1.5 in \cite{amgi13}]\label{teo1.5amgi}
	If the cost function $c\colon X\times Y$ is lower semi-continuous and bounded from below,
	then there exists a minimizer for \eqref{kantorovich} for any $\mu_0\in \mc P(X)$ and $\mu_1\in\mc P(Y)$.
\end{teo}

\begin{rmk}
	If $X=Y$ and $c(x,y)=d^2(x,y)$, the Kantorovich problem induces a distance $W_2(\mu,\nu)=\sqrt{C_{\Kan}(\mu,\nu)}$
	called Wasserstein distance, and $(\mc P_2(M),W_2)$ is
	a complete metric space, where $\mc P_2(M)$ is the space of probabilities
	with finite $2$-momentum.
	
	In this manuscript we will usually work with $X=Y=M$ a sub-Riemannian manifold
	and $c(x,y)=d^2_{SR}(x,y)$ the square distance. The distance is continuous,
	therefore the conditions for the existence of a minimizer are attained.\newline
\end{rmk}

%
%

Up to this point, we just used the metric structure on $M$.
If we look also at the differential structure on $M$, we can get a dynamical formulation
of the Optimal Transport problem. If $M$ is any complete Riemannian manifold,
we introduce the Continuity Equation, describing the evolution of a measure path in $\mc P(M)$
along the flow of a family of vector fields. If $I=[0,1]$,
consider a Borel family of probability measures
$\mu\colon I\to \mc P(M)$ and a Borel family of vector fields $v\colon I\times M\to TM$
respecting the condition
\begin{equation}\label{cond:v}
	\int_I\int_M\norm {v_t(x)}\d\mu_t(x)\d t<\infty.
\end{equation}
Then, the Continuity Equation states
\begin{equation}\label{eq:ceoriginal}
	\frac{\d}{\d t}\mu_t+\Div(v_t\mu_t)=0.
\end{equation}
This has to be considered in the distributional sense, meaning that
\begin{equation}\label{eq:ce}
	\int_I\int_M\de_t\phi(t,x)+\langle v_t,\nabla\phi(t,x)\rangle_x\,\d\mu_t(x)\d t=0\ \ \ \forall \phi\in C^\infty_c((0,1)\times M).
\end{equation}
This has been studied from many point of views,
see in particular \cite[Chapter 8]{ags08} and \cite{amcri14} for a wide overwiev.
Here, we are interested to one particular result, the so called \cor{superposition principle}.

Let $\Gamma:=C([0,1],M)$ be the set of continuous curves on $M$,
and $e_t\colon \Gamma\to M$ the map such that $\Gamma\ni\gamma\mapsto \gamma(t)$
for any $t\in I$.
If $(\mu_t,v_t)$ is a pair respecting the Continuity Equation~\eqref{eq:ceoriginal},
then  $\mu_t\in \mc P_2(M)$ is absolutely continuous for a.e.~$t\in I$ \cite[Theorem 2.29]{amgi13}
 (with respect to the Wasserstein distance)
and there exists $\widetilde \mu\in \mc P(\Gamma)$ concentrated on absolutely continuous curves
such that
\begin{equation}\label{superpo1}
	(e_t)_\#\widetilde \mu=\mu_t\quad \forall t\in I.
\end{equation}
We remark that in \cite{stetre17} this result have been extended to complete metric
spaces under opportune hypothesis.\newline

We consider now the adaptation of the Continuity Equation to the 
sub-Riemannian setting.
Given a sub-Riemannian variety $(M,H)$, we consider 
a Borel family of probability measures $\mu\colon I\to\mc P(M)$,
and a Borel vector field $v\colon I\times M\to HM$
such that for a.e.~$t\in I$, $v_t$ is an horizontal vector field on $M$.
Then, the Continuity Equation is obtained by just
replacing the gradient with the horizontal gradient,
\begin{equation}\label{eq:cesub}
	\int_I\int_M\de_t\phi(t,x)+\langle v_t,\nabla_H\phi(t,x)\rangle_x\,\d\mu_t(x)\d t=0\ \ \ \forall \phi\in C^\infty_c((0,1)\times M).
\end{equation}

\begin{rmk}\label{rmkce}
	For the Continuity Equation in Riemannian and
	sub-Riemannian setting, there are two equivalent 
	formulations that we state here without proof
	because they won't be needed in this work. 
	For a proof see \cite{ags08}.
	
	The ``weak'' formulation of the
	continuity equation considers smooth test functions 
	 defined over $M$,
	\begin{equation}\label{eq:cew}
		\frac{\d}{\d t}\int_M f\,\d\mu_t=-\int_M\langle v_t,\nabla f\rangle\, \d\mu_t\ \ \ \forall f\in C^\infty_c(M).
	\end{equation}
	If we take a function $\phi(t,x)=\theta(t)\cdot f(x)$ with $\theta\in C^\infty_c((0,1))$ and $f\in C^\infty_c(M)$,
	then  \eqref{eq:ce} implies  \eqref{eq:cew}  
	in the distributional sense.
	The two formulations are in fact equivalent, and furthermore
	\eqref{eq:cew} implies that
	it always exists a narrowly continuous family $\overline \mu\colon I\to \mc P(M)$,
	such that $\overline \mu_t=\mu_t$ for a.e.~$t\in [0,1]$
	and
	\begin{align}
		\mathclap{\int_{t_1}^{t_2}}\,\int_M\de_t\phi(t,x)+\langle v_t(x),\nabla_H\phi(t,x)\rangle \d\overline \mu_t(x)\d t&=\int_M\phi(t_2,x)\d\overline \mu_{t_2}(x)-\int_M\phi(t_1,x)\d\overline \mu_{t_1}(x)\label{eq:cepower}\\
		& \nonumber\\
		&\forall \phi\in C^\infty_c(I\times M),\ \forall t_1,t_2\in [0,1].\nonumber
	\end{align}
\end{rmk}

The Benamou--Brenier problem is a dynamical formulation of
the Optimal Transport problem. In the setting of a sub-Riemannian manifold
this becomes,
\begin{equation}\tag{BB}\label{BB}
	\begin{split}
		J_{\BB}(\mu_t,v_t)&:=\int_{I}\int_M\norm{v_t(x)}^2\d\mu_t(x)\d t\\
		\\
	C_{\BB}(\mu,\nu)&:=\inf_{(\mu_t,v_t)}J_{\BB}(\mu_t,v_t)
	\end{split}
\end{equation}
where the   infimum is taken among all the pairs of a Borel family of probability measures $\mu_t$,
and a Borel horizontal vector field $v_t$ respecting condition \eqref{cond:v}, such that
$(\mu_t,v_t)$ respects the sub-Riemannian
Continuity Equation \eqref{eq:cesub}.\newline

In the following we will prove that under the hypothesis that $M$ is a complete sub-Riemannian
with no boundary and no abnormal geodesics, if $\mu_0,\mu_1$ have finite $2$-momentum, then
the Kantorovich formulation and the Benamou--Brenier formulation are equivalent.
Namely, $C_{\Kan}(\mu_0,\mu_1)=C_{\BB}(\mu_0,\mu_1)$ and both
infima are realized.\newline

\section{The results}

In order to show that
the Kantorovich formulation \eqref{kantorovich} and the Benamou--Brenier
formulation \eqref{BB} of the transport problem are equivalent on sub-Riemannian manifolds
(under some assumptions that will be made precise later),
we  weaken the Benamou--Brenier problem by formulating
it over measures in $\mc P(I\times HM)$. We start by introducing the following set
for any Polish space~$(X,d)$,
\begin{equation}
	\Y(I;X):=\left\{\eta\in \mc P(I\times X):\ \pi_\#^{(t)}\eta=\mc L\right\}
\end{equation}
here $\pi^{(t)}\colon I\times X\to I$ is the projection on the first coordinate
while $\mc L$ is the Lebesgue measure on $I$.
We will denote by $\eta_t$ a disintegration of $\eta$ with respect to $\mc L$
for $t\in I$.

Any $\eta$ in $\Y(I;X)$ is called a \cor{Young measure} on $X$. There's a large
literature on this class of measures and when necessary we will  refer to \cite{berna08}.\newline

\subsection{The relaxed problem}\label{relax}
If $X=TM$ the tangent bundle to a smooth manifold $M$,
and $\pi\colon TM\to M$ the canonical projection,
then we call \cor{transport measure} a measure $\eta$~in~$\Y(I;TM)$
such that $\eta$ verifies the following version of the Continuity Equation,
\begin{equation}\label{transport}
\int_{I\times TM}\de_t\phi(t,\pi(v))+\langle v,\nabla \phi(t,\pi(v))_{\pi(v)}\,\d\eta(t,v)=0\quad \forall \phi\in C^\infty_c((0,1)\times M).
\end{equation}
\begin{rmk}
	Analogously to what happens with the standard Continuity Equation (see Remark \ref{rmkce}),
	in \cite[Lemma 5 and 6]{berna08} it is proved that if $\eta$ is a transport measure,
	then there exists a narrowly continuous family $\mu_t\in \mc P(M)$
	such that $\pi_\#\eta_t=\mu_t$ for a.e.~$t\in I$.
	In the following, $\mu_t$ will be this narrowly continuous curve
	of measures.

	Moreover, for any function $\phi\in C^1(I\times M)$
	bounded and Lipschitz and for any interval $[t_1,t_2]\subset I$,
	the following is respected
	\begin{equation}
		\int_{[t_1,t_2]\times TM}\de_t\phi+\langle v,\nabla_H\phi\rangle\, \d\eta=\int_M\phi(t_1,x)\d\mu_{t_1}(x)-\int_M\phi(t_2,x)\d\mu_{t_2}(x).
	\end{equation}
\end{rmk}

To relax the Benamou--Brenier formulation of the optimal transport,
we want to consider the case of transport measures that are
concentrated on the horizontal distribution, meaning that
$\eta\in \Y(I;HM)\subset \Y(I;TM)$. Observe that $\Y(I;HM)$ is closed
inside $\Y(I;TM)$ with respect to the narrow topology.

\begin{rmk}
	If $\eta$ is a  transport measure, then the 
	defining Continuity Equation can be formulated as
	\[
	\int_I\int_{TM}\left(\de_t\phi+v\phi\right)\,\d\eta(t,v)=0\quad\forall \phi \in C^\infty_c((0,1)\times M)
	\]
	and the sub-Riemannian setting appears only 
	with the additional condition that $\eta$ must be concentrated on $HM\subset TM$.
	Therefore, the results of \cite{berna08} may be used in our setting, without any
	further generalization, by considering the smooth structure on $M$.
\end{rmk}

\begin{rmk}
	We are working on a smooth manifold without boundary $M$ and we consider
	the tangent bundle $TM$. In the setting of Young measures, 
	$TM$ has to be endowed with a complete distance $D$ with the property
	that
	\[
	\frac{1+D((x_0,0),(x,v))}{1+||v||_x}\quad\m{with }v\in T_xM,
	\]
	is bounded on $TM$ for a given (and then every) point $x_0\in M$.
	The existence of such a distance is proved in \cite[\S3]{berna08}.
	Observe that we only need a topological compatibility with the sub-Riemannian structure
	on $M$, therefore we can consider the smooth structure
	on $M$ as specified in the previous remark.\newline
\end{rmk}

If we pose $\eta\in \Y(I,TM)$,
Equation \eqref{transport}
can be equivalently stated in the setting of smooth manifolds,
	\[
\int_I\int_M\de_t\phi+v\phi\,\d\eta(t,v)=0\quad\forall \phi \in C^\infty_c((0,1)\times M).
\]
If we set the additional condition that $\eta$ must be concentrated on $HM\subset TM$,
we have a transport measure in $\Y(I;HM)$. Therefore,
reversing the conditions,
we get a sub-Riemannian relaxed version of the Continuity Equation. Indeed,
a measure $\eta \in \Y(I;HM)$ 
is a transport measure if and only if
\begin{equation}\label{con-t}
	\int_{I\times HM}\de_t\phi(t,\pi(v))+\langle v,\nabla_H\phi(t,\pi(v))\rangle_{\pi(v)}\d\eta(t,v)=0\ \ \ \forall \phi\in C^\infty_c((0,1)\times M).
\end{equation} 
where $\nabla_H$ is the horizontal gradient.
A slightly stronger version links $\eta$ to its ``values'' at the extremes of the interval,
\begin{equation}\label{con-star}
	\int_{I\times HM}\de_t\phi+\langle v,\nabla_H\phi\rangle\d\eta=\int_M\phi(1,x)\d\mu_1(x)-\int_M\phi(0,x)\d\mu_0(x) \ \ \ \forall \phi\in C^\infty_c(I\times M).
\end{equation}
These versions of the Continuity Equation define two closed subsets of
$\Y(I;HM)$ with respect to the narrow topology,
\begin{align*}
	\mc T_H((0,1);M)&:=\left\{\eta\in\Y(I;HM):\ \eta\ \m{verifies \eqref{con-t}}\right\}\\
	\Y_{\mu_0}^{\mu_1}(I;HM)&:=\left\{\eta\in \Y(I;HM):\ \eta\ \m{verifies \eqref{con-star} w.r.t.}\ \mu_0,\mu_1\right\}.
\end{align*}
Observe that $\Y_{\mu_0}^{\mu_1}(I;HM)\subset\mc T_H((0,1);M)$.
We can now formulate
the relaxed version of~\eqref{BB},
\begin{equation}\tag{BB$^\star$}\label{BBstar}
	\begin{split}
		J_{\BB}^\star(\eta)&:=\int_{I\times HM}\norm{v}^2\,\d\eta(t,v)\\
		\\
		C_{\BB}^\star(\mu_0,\mu_1)&:=\inf_{\eta\in\Y^{\mu_1}_{\mu_0}}J_{\BB}^\star(\eta).
	\end{split}
\end{equation}

\begin{rmk}\label{weakening}	
We point out that the relaxed version is in fact a \cor{weaker} version of \eqref{BB}
because to any pair $(\mu_t,v_t)$ satisfying the Continuity Equation \eqref{eq:ce},
we can naturally associate the measure $\eta_{(\mu,v)}=\mu_t\xx\delta_{v_t}\in \Y_{\mu_0}^{\mu_1}(I;HM)$,
that is the measure defined by
\[
\phi\mapsto \int_{I}\int_M\phi(t,v_t(x))\, \d\mu_t(x)\d t\ \ \ \forall \phi\in C_c^\infty(I\times HM).
\]
We have $J_{\BB}^\star(\eta_{(\mu,v)})=J_{\BB}(\mu_t,v_t)$ and therefore
we are in fact weakening the constaints on the \eqref{BB} problem.\newline
\end{rmk}

Recall that $\pi\colon HM\to M$ is the natural projection from the horizontal bundle
to the base manifold.
By Lemma \ref{lem:etamu}, 
if $\eta\in \mc T_H((0,1);M)$ is a transport measure, then there exists
a narrowly continuous family $\mu_t^\eta\colon I\to \mc P(M)$ 
and a disintegration $\eta_t$ of $\eta$ such that 
\begin{equation}\label{eq:etamu}
	\pi_\#\eta_t=\mu_t^\eta\quad \forall t\in I.
\end{equation}

Furthermore, for any $\eta\in \mc T_H((0,1);M)$ and
 for any $t\in I$,
we consider a disintegration $\eta_{t,x}$ with respect to $\mu_t^\eta$.
Through $\eta_{t,x}$ we can introduce
an associated measurable map 
\[
v^\eta\colon I\times M\to HM
\]
 such that 
for a.e.~$t$, $ v_t^\eta$ is a horizontal vector field on $M$,
\begin{equation}\label{eq:etav}
	 v_t^\eta(x):=\int_{H_xM}v\, \d\eta_{t,x}(v).
\end{equation}
Observe that by construction we get
\begin{equation}\label{equal}
	\begin{split}
	\int_{I\times HM}\left(\de_t\phi+\langle v, \nabla_H\phi\rangle\right)\d \eta(t,v)
	&=\int_{I\times M}\int_{H_xM}\left(\de_t\phi+\langle v,\nabla_H\phi\rangle \right)\d\eta_{t,x}(v)\d\mu_t^\eta(x)\d t\\
	&=\int_{I\times M}\left(\de_t\phi+\langle v^\eta_t,\nabla_H\phi\rangle\right)\d \mu^\eta_t(x)\d t\\
	&\ \forall \phi\in C^\infty_c(I\times M)
	\end{split}
\end{equation}
where we used the definition of $ v^\eta_t$ and the linearity of the scalar product.

Moreover, if $J^\star_{\BB}(\eta)<+\infty$, then
by Jensen's inequality
\[
+\infty>\int_{I\times HM}\|v\|^2\d\eta\geq \int_I\int_M\|v_t^\eta\|^2\d\mu^\eta_t\d t
\]
that is, $(\mu^\eta_t,v^\eta_t)$ verifies condition \eqref{cond:v}. \newline

Therefore, we get two results that we summarize in the following remark.
\begin{rmk}\label{rmkmuv}
From \eqref{equal} and the fact that $\eta$ satisfies \eqref{con-t},
we get that $(\mu_t^\eta,v_t^\eta)$ satisfies the Continuity Equation \eqref{eq:cesub}.
Also, observe that if we impose the narrow continuity of the curve $\mu_t^\eta$,
then it is uniquely defined and 
$\eta\in \Y_{\mu_0}^{\mu_1}(I;HM)$, where
we used the notation $\mu_0,\mu_1$ for $\mu_0^\eta,\mu_1^\eta$ on the indices,
for the sake of clarity.

To any Young transport measure $\eta$ concentrated on the horizontal bundle, we have thus associated
a unique pair $(\mu_t,v_t)$ of a narrowly continuous curve of measures
and a Borel family of vector fields (respecting \eqref{cond:v})
that satisfy the Benamou--Brenier constraint \eqref{eq:cesub}.\newline
\end{rmk}

\begin{teo}\label{existence}
	Consider $\mu_0,\mu_1\in \mc P_2(M)$,
	suppose that there exists $\eta\in \Y_{\mu_0}^{\mu_1}(I;HM)$ with
	$J_{\BB}^\star(\eta)$ finite (\cor{i.e.}~the relaxed problem is feasible).
	Then, 
	\[
	C_{\BB}(\mu_0,\mu_1)=C_{\BB}^\star(\mu_0,\mu_1).
	\]
\end{teo}
\proof 
As a consequence of Remark \ref{weakening}, we know that
\[
C_{\BB}(\mu_0,\mu_1)\geq C_{\BB}^\star(\mu_0,\mu_1). 
\]
We are going to prove the reverse inequality,
in particular for any transport measure $\eta$
we show that that the pair $(\mu^\eta_t,v^\eta_t)$
has a $J_{\BB}$-cost smaller than $J_{\BB}^\star(\eta)$.\newline

Let $\eta\in \Y_{\mu_0}^{\mu_1}(I;HM)$,
we consider
$\mu_t^\eta$ and $v_t^\eta$ introduced in \eqref{eq:etamu} and \eqref{eq:etav} respectively.
As summarized in Remark \ref{rmkmuv},
the pair $(\mu_t^\eta,v_t^\eta)$ satisfies the continuity equation \eqref{eq:ce}.
Moreover, the map $v\mapsto\norm{v_t}_x^2$ is convex for any $(t,x)\in I\times M$,
	therefore by Jensen's inequality,
	\begin{align*}
		J_{\BB}(\mu_t^\eta,v_t^\eta)=\int_{I\times M}\norm{v_t^\eta(x)}^2\,\d\mu^\eta_t(x)\d t
		&\leq \int_{I\times M}\int_{H_xM}\norm{v}^2_x\d\eta_{t,x}(v)\,\d\mu^\eta_t(x)\d t\\
		&=\int_{I\times HM}\norm{v}^2\,\d\eta(t,v)\\
		&=J_{\BB}^\star(\eta)
	\end{align*} 
	This proves $C_{\BB}(\mu_0,\mu_1)\leq C_{\BB}^\star(\mu_0,\mu_1)$ and 
	allows to complet the proof.\fine
	
	\subsection{Equivalence with the Kantorovich formulation}\label{kanequiv}
	It remains to prove that the Kantorovich formulation \eqref{kantorovich}
	and the Benamou--Brenier  formulation \eqref{BB} are equivalent.
	
	We denote by $\Geod(M)$ the set of constant speed geodesics on the sub-Riemannian manifold $(M,H)$,
	with the topology induced by the sup distance 
	$\sup_{t\in[0,1]}d(\omega^{(1)}(t),\omega^{(2)}(t))$
	between any two geodesics $\omega^{(1)},\omega^{(2)}$.\newline
	
	The key technical result that we are going to use is that for any Borel measure $\rho$ on $M\times M$,
	there exists a Borel map 
	\[
	S\colon M\times M\to \Geod(M)
	\]
	defined $\rho$-a.e.~such that 
	\[
	S_0(x,y)=x\ \ \ \m{and}\ \ \ S_1(x,y)=y.
	\]
	This result is proved in Theorem \ref{teoSuslin} in the Appendix.\newline
	\color{black}
	
	The map $S$ above allows to build a map $\Pi(\mu_0,\mu_1)\to \Y_{\mu_0}^{\mu_1}(I;HM)$.
	Indeed, consider the measurable function
	$F\colon I\times M\times M\to I\times HM$
	such that
	\begin{equation}\label{Fmap}
	(t,x,y)\mapsto (t,\de_tS(x,y)).
	\end{equation}

	Then, for any admissible plan $\gamma\in \Pi(\mu_0,\mu_1)$, we consider the map
	\begin{equation}\label{mapF}
		\gamma\mapsto F_\#(\mc L\xx\gamma).
	\end{equation}
	\begin{lemma}
		For any $\gamma\in \Pi(\mu_0,\mu_1)$, the image of the map
		\eqref{mapF} is in $ \Y^{\mu_1}_{\mu_0}(I;HM)$.
		\end{lemma}
		\proof
		Throughout this proof we will use the notation $\eta:=F_\#(\mc L\xx \gamma)$.
		First observe that $\pi^{(t)}\circ F=\id_I$, therefore $\pi^{(t)}_\#\eta=\mc L$
		and so $\eta\in \Y(I;HM)$.
		
		For any test function $\phi\in C^\infty_c(I\times M)$ we evaluate the following derivative for any 
		pair of points $(x,y)\in M\times M$,
		\begin{equation}\label{eqdt}
			\frac{\d}{\d t}\phi(t,S_t(x,y))
			= \de_t\phi(t,S_t(x,y))+\langle \de_tS_t(x,y),\nabla_H\phi(t,S_t(x,y))\rangle,
		\end{equation}
		where we used the fact that $S_t$ is admissible and therefore
		$\de_tS_t$ is almost always an horizontal vector in $H_{S_t}M$.
		Therefore, by construction we get
		\begin{align*}
			\int_{I\times HM}\de_t\phi(t,\pi(v))&+\langle v,\nabla_H\phi(t,\pi(v))\rangle\,\d\eta(t,v)\\
			&=\int_{I\times M\times M}\de_t\phi(t,S_t(x,y))+\langle \de_tS_t(x,y),\nabla_H\phi(t,S_t(x,y))\rangle\d\gamma(x,y)\d t\\
			&\m{using \eqref{eqdt}}\\
			&=\int_{M\times M}\phi(1,x)-\phi(0,y)\d\gamma(x,y)\\
			&=\int_M\phi(1,x)\d\mu_1(x)-\int_M\phi(0,x)\d\mu_0(x).
		\end{align*}
		This is precisely the condition we wanted to prove.\fine
		
	\begin{lemma}\label{lemcsg}
		If $F$ is defined as above and $\gamma\in \Pi(\mu_0,\mu_1)$, then
	\[
		J_{\Kan}(\gamma)=J_{\BB}^\star(F_\#(\mc L\xx\gamma)).
	\]
	\end{lemma}
	\proof
	Recall that the map $S$ associates to any pair of points $(x,y)$
	a constant speed geodesic between them. Therefore
	\[
	d^2(x,y)=\int_I\norm{\de_tS_t(x,y)}^2\d t.
	\]
	As $\gamma$ is an admissible plan between $\mu_0$ and $\mu_1$ we pass to
	the evaluation of the Kantorovich cost of $\gamma$,
	\begin{align*}
		J_{\Kan}(\gamma)=\int_{M\times M}d^2(x,y)\d \gamma
		&=\int_{I\times M\times M}\norm{\de_tS_t(x,y)}^2\d\gamma\d t\\
		&=\int_{I\times HM}\norm{v}^2\d\eta(t,v)=J_{\BB}^\star(\eta),
	\end{align*}
	where again we used the notation $\eta=F_\#(\mc L\xx \eta)$.\fine
	
	A direct corollary of this lemma is the inequality
	\begin{equation}\label{starlesskan}
		C_{\BB}^\star(\mu_0,\mu_1)\leq C_{\Kan}(\mu_0,\mu_1).
	\end{equation}
	In particular, observe that the feasibility of Kantorovich problem
	implies the feasibility of the relaxed Benamou--Brenier problem.\newline
	
	In order to prove the last main result,
	we develop in the Appendix section \ref{superyoung}
	the theory of \cor{generalized curves} and
	of transport measure decompositions.
	
	\begin{defin}\label{def_gc}
		A generalized curve is a transport
		measure $\eta$ with $\mu=\pi_\#\eta$ such that 
		$\mu_t$ is a Dirac measure for a.e.~$t\in I$.
		 In particular, there exists $\omega\colon I\to M$
		continuous such that $\mu_t=\delta_{\omega(t)}$. The set of generalized curves
		is denoted by $\mc G(I,M)$ and is closed inside $\mc T((0,1);TM)$.
	\end{defin}
	
	As made precise in the Appendix,
	any such curve $\omega$ is absoluteluy continuous,
	and we denote~by 
		\[
	\m{pr}\colon \mc G(I,M)\to AC(I,M)
	\]
	the projection from the space of generalized curves to
	absolutely
	continuous curves on $M$.

	\begin{teo}\label{kanstar}
		Given  $\mu_0,\mu_1\in \mc P_2(M)$, 
		the Kantorovich formulation and
		the Benamou--Brenier formulation are equivalent,
		\[
		C_{\Kan}(\mu_0,\mu_1)=C_{\BB}^\star(\mu_0,\mu_1)=C_{\BB}(\mu_0,\mu_1)
		\]
	\end{teo} 
	\proof
	After Theorem \ref{existence} and the inequality \eqref{starlesskan}, it remains to prove
	that 
	\[
	C_{\Kan}(\mu_0,\mu_1)\leq C_{\BB}(\mu_0,\mu_1)=C_{\BB}^\star(\mu_0,\mu_1).
	\]
	Let $e_t\colon C(I,M)\to M$ be the map that sends any continuous curve $\gamma\colon I\to M$ to the point~$\gamma(t)$.
	We consider the map
	\begin{equation}\label{Emap}
	E:=(e_0,e_1)\circ \m{pr}\colon \mc G_H(I,M)\to M\times M.
	\end{equation}
	Here $\mc G_H(I,M)=\m{pr}^{-1}(\Omega_{HM})$ the
	space of generalized curves above admissible curves with respect to the
	sub-Riemannian structure.
	
	The Young superposition principle \ref{superpoyoung},
	states that for every $\eta\in \mc T_H((0,1);M)$
	there exists a decomposition measure $\widetilde \eta\in\mc P(\mc G_H(I,M))$
	such that 
	\[
	E_\#\widetilde \eta=\pi_\#\eta.
	\]
	Observe in particular that if $\eta\in \mc Y^{\mu_1}_{\mu_0}(I;HM)$,
	then 
	 $E_\#\widetilde \eta\in \Pi(\mu_0,\mu_1)$.
	Also, we use the notation
     $\omega=\m{pr}(\nu)$ for the curve associated to a generalized curve~$\nu$. Then,
	\begin{align*}
		\int_{M\times M}d^2(x,y)\d E_\#\widetilde \eta
		&=\int_{\mc G_H}d^2(\omega_0,\omega_1)\d\widetilde \eta(\nu)\\
		&\leq\int_{\mc G_H}\int_I\norm{\dot\omega_t}^2\d t\d\widetilde \eta(\nu)\\
		&\m{(because  of \eqref{vtnut} and convexity)}\\
		&\leq\int_{\mc G_H}\int_I\int_{H_{\omega_t}M}\norm{v}^2\d\nu_t\d t\d\widetilde \eta(\nu)\\
		&\m{(because of \eqref{etadec} applied on $\widetilde\eta$)}\\
		&=\int_{I\times HM}\norm{v}^2 \d\eta(t,v)=J_{\BB}^\star(\eta).
	\end{align*}
	This proves that $J_{\Kan}(E_\#\widetilde \eta)\leq J_{\BB}^\star(\eta)$
	and therefore in particular
	\[
	C_{\Kan}(\mu_0,\mu_1)\leq C_{\BB}^\star(\mu_0,\mu_1)=C_{\BB}(\mu_0,\mu_1).
	\]\fine

Resuming our results, we proved in Theorem \ref{existence}
that $C_{\BB}(\mu_0,\mu_1)=C_{\BB}^\star(\mu_0\mu_1)$,
and in Theorem \ref{kanstar} that $C_{\Kan}(\mu_0,\mu_1)=C_{\BB}^\star(\mu_0,\mu_1)$.
Moreover, thanks to the latter theorem and Lemma \ref{lemcsg}
we know that if $\gamma$ is an optimal Kantorovich plan, then $\eta=F_\#{(\mc L\xx\gamma)}$ 
realizes the optimal Benamou--Brenier cost.

In order to complete the proof of Theorem \ref{mainteo}
it remains to show that if $(\mu_t,v_t)$ is a minimizer for Benamou--Brenier
and $\widetilde\eta\in \mc P(\mc G_H(I,M))$ a decomposition in the sense of \eqref{etadec},
then $\widetilde \eta$ is supported on geodesics.

\begin{teo}
	Suppose the $(\mu_t,v_t)$ respects the condition \eqref{cond:v}
	and the Continuity Equation \eqref{eq:ceoriginal},
	and moreover
	\[
	J_{\BB}(\mu_t,v_t)=C_{\BB}(\mu_0,\mu_1),
	\]
	meaning that $(\mu_t,v_t)$ is a minimizer for Benamou--Brenier.
	Consider $\widetilde\eta\in \mc P(\mc G_H(I,M))$ a decomposition for $\mu_t$,
	then $\m{pr}(\nu)$ is a geodesic for $\widetilde\eta$-a.e.~generalized curve $\nu$.
\end{teo}

\proof We suppose there exists a measurable subset $A\subset \mc G_H(I,M)$
such that $\widetilde \eta(A)>0$ and such that 
\[
\m{pr}(A)\cap \Geod(M)=\varnothing.
\]
Then, 
we are going to
build a map $\mc G_H(I,M)\to \mc G_H(I,M)$
that sends the decomposition $\widetilde \eta$
to a measure $\widetilde \eta^*$ inducing 
a probability $\eta^*$ with a smaller Benamou--Brenier cost than $\eta$.\newline

We consider the Borel map $S\colon M\times M\to \Geod(M)$ from Theorem \ref{teoSuslin}
and we recall that it is defined $(\mu_0\xx\mu_1)$-a.e.
From this we can
build the map $F\colon M\times M\to \mc G_H(I,M)$
such that
\[
(x,y)\mapsto \d t\xx\delta_{S_t(x,y)}\xx \delta_{\de_tS_t(x,y)},
\]
that is the generalized curve concentrated on the tangent vectors to the (sub-Riemannian)
geodesic between $x$ and $y$. Observe that if $\nu\in A$ and $\omega=\m{pr}(\nu)$
is an admissible curve that is not a (constant speed) geodesic,
then
\begin{equation}\label{eqlength}
	\int_I \norm{\dot\omega_t}^2\d t\geq\ell(\omega)^2>\ell(S(\omega_0,\omega_1))^2= \int_I\norm{\de_t S_t(\omega_0,\omega_1)}^2\d t,
\end{equation}
where we used the minimization property of geodesics.\newline

We introduce the measurable map $m\colon \mc G_H(I,M)\to \mc G_H(I,M)$
such that
\[
\left\{
\begin{array}{ll}
	m|_A&= F\circ (e_0,e_1)\circ\m{pr}\\
	m|_{\mc G_H\backslash A}&=\id.
\end{array}
\right.
\]
We consider the modified measure $\widetilde \eta^*:=m_\#\widetilde \eta$
and also the probability
\[
\eta^*:=\int_{\mc G}\nu\,\d\widetilde \eta^*(\nu).
\]
By construction $\eta^*\in \Y_{\mu_0}^{\mu_1}(I;HM)$
and $\widetilde \eta^*$ is a decomposition of $\eta^*$.
Using the defintion \eqref{etadec} of decomposition, we develop
the following calculations where again we use the notation $\omega=\m{pr}(\nu)$,
meaning that $\nu$ is a generalized curve above the continuous curve $\omega$.
\begin{align*}
	J^\star_{\BB}(\mu_t\xx \delta_{v(t)})=\int_{I\times M}\norm{v_t}^2\d \mu_t \d t
	&=\int_{\mc G_H}\int_{I\times HM}\norm{v}^2\d\nu(v)\d\widetilde \eta(\nu)\\
	&=\int_{\mc G_H\backslash A}\int_{I\times HM}\norm{v}^2\d\nu\d\widetilde \eta
	+\int_A\int_{I\times HM}\norm{v}^2\d\nu\d\widetilde \eta\\
	&\m{(by \eqref{eqnu})}\\
	&=\int_{\mc G_H\backslash A}\int_{I\times HM}\norm{v}^2\d\nu\d\widetilde \eta
	+\int_A\int_{I\times HM}\norm{v}^2\d t \d \nu_t\d \widetilde \eta\\
	&\m{(by convexity)}\\
	&\geq\int_{\mc G_H\backslash A}\int_{I\times HM}\norm{v}^2\d\nu\d\widetilde \eta
	+\int_A\int_I{\left|\left|\int_{H_{\omega_t}M}v\,\d\nu_t\right|\right|}^2\d t\d\widetilde \eta\\
	&\m{(by \eqref{vtnut})}\\
	&=\int_{\mc G_H\backslash A}\int_{I\times HM}\norm{v}^2\d\nu\d\widetilde \eta
	+\int_A \int_I\norm{\dot\omega_t}^2\d t\d\widetilde \eta\\
	&\m{(because of \eqref{eqlength})}\\
	&>\int_{\mc G_H\backslash A}\int_{I\times HM}\norm{v}^2\d\nu\d\widetilde \eta
	+\int_A \int_I\norm{\de_tS_t(\omega_0,\omega_1)}^2\d t \d\widetilde \eta\\
	&=\int_{\mc G_H\backslash A}\int_{I\times HM}\norm{v}^2\d\nu\d\widetilde \eta
	+\int_{\overline A} \int_{I\times HM}\norm{v}^2\d\nu\d\widetilde \eta^*\\
	&=\int_{\mc G_H}\int_{I\times HM}\norm{v}^2\d\nu\d\widetilde \eta^*=J^\star_{\BB}(\eta^*)
\end{align*}
where we used the notation $\overline A$ for the image $\overline A:=F\circ (e_0,e_1)\circ \m{pr}(A)\subset \mc G_H(I,M)$.

Therefore, $J_{\BB}(\mu_t,v_t)>J_{\BB}^\star(\eta^*)$
which is an absurd because $(\mu_t,v_t)$ was
supposed to be a minimizer for the Benamou--Brenier formulation.
This proves that $\widetilde \eta$ must be supported on geodesics.\fine

	\appendix
	\section*{Appendix}\label{Some technical results}
	\addcontentsline{toc}{section}{Appendix}
	\addtocounter{section}{1}
	\setcounter{teo}{0}
	\setcounter{equation}{0}
	\subsection{Transport measures}\label{transpmeas}

\begin{lemma}\label{lem:etamu}
	Consider $\eta\in \mc T_H((0,1);M)$  a transport measure over $(M,H)$, then
	there exists a narrowly continuous family $\mu_t\colon I\to\mc P(M)$
	and a disintegration $\eta_t$ of $\eta$, such that
	$\pi_{\#}\eta_t=\mu_t$ for each $t\in I$.
	As a consequence $\eta\in\Y_{\mu_0}^{\mu_1}(I;HM)$.
\end{lemma}
\proof
Consider $\widetilde \eta_t$  a disintegration of~$\eta$,
and $\widetilde \mu_t=\pi_\#\widetilde \eta_t$.
Then, by \cite[Lemma 37]{berna08} $\widetilde \mu_t$ is equal for a.e.~$t$
to a narrowly continuous map $\mu_t$
if and only if for any compactly supported smooth function $f\colon I\to \R$,
the following function~$I\to \R$,
\[
F(t):=\int_M f\,\d\widetilde \mu_t,
\]
is equal a.e.~to a continuous function.

Consider Equation \eqref{con-t} and apply it to a function $\phi(t,x):=\theta(t)\cdot f(x)$
where $\theta\in C^\infty_c((0,1))$. We obtain
\[
\begin{split}
0&=\int_I\int_{HM} \theta'(t)f(\pi(v))+\theta(t)\langle v, \nabla_H f(t,\pi(v))\rangle\d\widetilde \eta_t(v)\d t\\
&=\int_I\int_M \theta'(t)f(x)\d\widetilde \mu_t(x)+\int_I\theta(t)\left(\int_{HM}\langle v,\nabla_Hf(\pi(v))\rangle \d\widetilde \eta_t(v)\right)\d t\\
&=\int_I\theta'(t)F(t)+\theta(t)\left(\int_{HM}\langle v,\nabla_Hf\rangle\d\widetilde \eta_t\right)\d t.
\end{split}
\]
As the choice of $\theta$ is arbitrary, this means that $F'(t)=\int_{HM}\langle v,\nabla_H f\rangle \d\widetilde \eta$ in the sense of distributions,
therefore $F$ is equal a.e.~to an absolutely continuous function.\fine

\subsection{The $S$ map}\label{Smap}

In order to 
treat the superposition principle in sub-Riemannian setting, we
need to prove the existence of a measurable map that to any pair of points on $M$
associates a geodesics connecting the two points. For this, we need
to use some results from the theory of Suslin spaces,  we will use
\cite{boga07} as the main reference on this subject.

\begin{defin}
	A set in a Hausdorff topological space is a Suslin set
	if it is the image of Polish space 
	with respect to a continuous map.
	
	A Suslin space is a Hausdorff topological space that is also a Suslin set.
\end{defin}

On a Hausdorff space $X$, we use the notation $\mc B(X)$
for its Borel $\sigma$-algebra.
In $X$, every Borel set is also Suslin.
The converse is not true, but there exists a notion similar to measurability
that is respected by any Suslin set. 

Given a Borel measure $\rho$ on $X$, a set function $\rho^*$ 
is defined for any subset $A\subset X$,
\[
\rho^*(A):=\inf\left\{\sum_{i=1}^\infty\rho(A_i):\ A_i\m{ is Borel},\ A\subset \bigcup_{i=1}^\infty A_i\right\}.
\]
Observe that this notion can be introduced
for any set function $\rho$ not necessarily a Borel measure, anyway
we will use this definition with the aim of simplifying the subject.

\begin{defin}
	A subset $A\subset X$ is $\rho$-measurable or measurable with respect to $\rho$,
	if for any $\varepsilon >0$ there exists a Borel set $A_\varepsilon\in \mc B(X)$
	such that
	\[
	\rho^*(A\triangle A_\varepsilon)<\varepsilon.
	\]
\end{defin}

\begin{prop}[See {\cite[Corollary 1.5.8]{boga07}}]\label{propsousmeas}
	Consider a Borel measure $\rho$ on $X$.
	A set $A\subset X$ is $\rho$-measurable if and only if
	there exists $A^{(1)},A^{(2)}\in \mc B(X)$ 
	such that
	\[
	A^{(1)}\subset A\subset A^{(2)}\quad\m{and}\quad \rho(A^{(2)}\backslash A^{(1)})=0.
	\]
\end{prop}

Now we can state a general property of Suslin sets.
\begin{teo}[See {\cite[Theorem 7.4.1]{boga07}}]\label{teosousmeas}
	If $X$ is a Hausdorff space and $\rho$ a Borel measure on~it,
	then every Suslin set $A\subset X$ is $\rho$-measurable.\newline
\end{teo}

Here we can state the main result of this appendix section.
\begin{teo}\label{teoSuslin}
	Consider a sub-Riemannian complete connected manifold $M$
	without boundary and with no non-trivial abnormal geodesics, and
	the (Hausdorff) space of constant speed geodesics on it $\Geod(M)$
	endowed with the sup distance.
	For any Borel measure $\rho$ on $M\times M$
	there exists a measurable map
	\[
	S\colon M\times M\to \Geod(M)
	\]
	defined $\rho$-a.e.~and such that 
	\[
	S_0(x,y)=x,\quad S_1(x,y)=y\quad \m{for }\rho-{a.e.}\,(x,y)\in M\times M.
	\]
\end{teo}
\proof
We introduce the map $\Psi$ sending every pair of points $(x,y)$
to the set of constant speed geodesics from $x$ to $y$.
We are going to prove that the hypothesis of \cite[Theorem 6.9.2]{boga07}
are verified, in particular the graph of $\Psi$, seen as a subset
of $M\times M\times AC([0,1],M)$, is a Polish space,
where we put the $L^\infty$ topology on the space $AC$.

Moreover,
we recall that any normal extremal in $AC([0,1],M)$ is uniquely determined
by its starting point $x\in M$ and a covector $\lambda \in T_x^*M$.
Indeed, normal geodesics are solutions
to a Hamiltonian system in $T^*M$ with some covector $\lambda$ as initial point.
Therefore, if $M$ has no non-trivial abnormal geodesics, then the map $T^*M\to AC([0,1],M)$ is continuous and 
has the subset of all the Pontryagin extremal trajectories as image. We observe
that for any covector $\lambda\in T_x^*M$ such that $H(\lambda)=0$, then the map send $\lambda$
to the constant geodesic $\gamma(t)\equiv x$. We then just proved that the subset
\[
\mc E:=\left\{(x,y,\gamma):\ \gamma\in AC([0,1],M)\m{ extremal tajectory},\,\gamma(0)=x,\,\gamma(1)=y\right\},
\]
is a Suslin set, because it is the image of a continuous mapping in a Polish space.

Moreover, there exists a continuous map
\[
L\colon \mc E\to M\times M\times \R
\]
sending any triplet $(x,y,\gamma)$ to $(x,y,\ell(\gamma))$.
Therefore the $\Psi$ graph is the preimage via $L$
of the closed set $\{(x,y,d(x,y))\}$,
thus the graph is also a Suslin set.\newline

As a consequence of \cite[Theorem 6.9.2]{boga07}, a selection $S\colon M\times M\to \Geod(M)$
exists, such that
\[
S(x,y)\in \Psi(x,y)\quad \forall (x,y)\in M\times M
\]
and $S$ is measurable with respect to the $\sigma$-algebra $\mc B(\Geod(M))$ of Borel sets
on $\Geod(M)$, and the $\sigma$-algebra $\sigma(S_{M\times M})$ generated by the Suslin
sets on $M\times M$.\newline

We conclude by showing that for any Borel measure $\rho$ on $M\times M$,
there exists a Borel set $B\subset M\times M$ with $\rho(B)=0$ such that the
restriction 
 $\left.S\right|_{M\times M\backslash B}$ is
(Borel) measurable.

Given a countable base $A_1,A_2,\dots$ of the Borel $\sigma$-algebra $\mc B(\Geod(M))$,
by Theorem \ref{teosousmeas} and Proposition \ref{propsousmeas}
for every $i=1,2,\dots$ there exist $B_i^{(1)},B_i^{(2)}$ both Borel sets in $M\times M$
such that
\[
B_i^{(1)}\subset S^{-1}(A_i)\subset B_i^{(2)}\quad\m{and}\quad \rho(B_i^{(2)}\backslash B_i^{(1)})=0.
\]
If we take 
\[
B:=\bigcup_{i=1}^\infty (B_i^{(2)}\backslash B_i^{(1)}),
\]
then $\rho(B)=0$ and $S$ is (Borel) measurable if restricted to $M\times M\backslash B$.\fine
\color{black}

\subsection{Young sub-Riemannian superposition}\label{superyoung}

We introduced generalized curves with Definition \ref{def_gc}.
From \cite[Lemma 8]{berna08} we get that for any generalized curve $\nu$, there exists an absolutely
continuous curve $\omega$ such that $\nu$ is a curve above $\omega$
and there exists a Borel family $\nu_t$ of probabilities on $T_{\omega_t}M$ for a.e.~$t\in I$,
such that 
\begin{equation}\label{eqnu}
\nu=\d t\xx \delta_{\omega_t}\xx\nu_t.
\end{equation}
Moreover, it is proved
that 
\begin{equation}\label{vtnut}
\int_{T_{\omega_t}M}v\,\d\nu_t(v)=\dot\omega_t\quad\m{for a.e.}\ t\in I.
\end{equation}
Viceversa, in order for \eqref{eqnu} to define a generalized
curve above $\omega$, it is necessary and sufficient that
\begin{itemize}
	\item $t\mapsto \int_{T_{\omega_t}M}\norm{v}\,\d \nu_t$ is Lebesgue integrable on $I$
	\item  \eqref{vtnut} is verified.\newline
\end{itemize}

We
consider the projection 
\[
\m{pr}\colon\mc G(I,M)\to AC(I,M)
\]
that to any generalized curve above $\omega$ associates the curve $\omega$.
By \cite[Theorem 9]{berna08}, this map is continuous. 

We denote by $\mc G_H(I,M)$ the set of horizontal generalized curves,
that is $\m{pr}^{-1}\left(\Omega_{HM}\right)$, meaning that any horizontal
generalized curve
$\nu$ is above an absolutely
continuous curve $\omega$ such that $\dot\omega_t$ lies
in the horizontal distribution $H_{\omega_t}M$ for a.e.~$t\in I$.

If $\eta\in \Y(I,TM)$ is a transport measure,
we know from \cite[Theorem 19]{berna08}, that there exists unique (up to a negligible set)
a Borel probability measure $\widetilde \eta$  on $\mc G(I,M)$ such that 
\begin{equation}\label{etadec}
\eta=\int_{\mc G}\nu\,\d\widetilde \eta(\nu).
\end{equation}
The measure $\widetilde \eta$ is called a \cor{decomposition} of $\eta$.\newline

Considering the results in the same work by Bernard, by \cite[Proposition 20]{berna08}
we have that for every transport measure $\eta$ concentrated on $I\times HM\subset I\times TM$,
if $\widetilde \eta$ is a decomposition then $\widetilde\eta$-a.e.~generalized curve is concentrated on $I\times HM$.
By \eqref{vtnut} this means that $\widetilde\eta$-almost every~generalized curve is horizontal, or equivalently
that $\widetilde \eta$ is concentrated on $\mc G_H(I,M)$. Indeed, if $\nu\in \mc G(I,M)$ is a generalized curve
above $\omega\in AC(I,M)$, and  is concentrated
on $I\times HM$, this means that
\[
\dot\omega_t=\int_{T_{\omega_t}M}v\,\d\nu_t\in H_{\omega_t}M\quad \m{for a.e.}\ t\in I.
\]

We resume this result and \cite[Proposition 21]{berna08} from the same article by Bernard, in the following statement.

\begin{teo}[Young superposition principle in sub-Riemannian setting]\label{superpoyoung}
	If $\eta$ is a transport measure concentrated on the horizontal distribution,
	\cor{i.e.}~$\eta\in \mc T_H((0,1);M)$, then a decomposition
	$\widetilde\eta$ of $\eta$ is concentrated on the set of (generalized) admissible curves $\mc G_H(I,M)$~and 
	\[
	(e_t\circ \m{pr})_\#\widetilde \eta=\mu_t=\pi_\#\eta_t.
	\]\newline
\end{teo}

\subsection{Compactness result}\label{cptresult}
Finally we prove a compactness result for 
the sublevels of the functional $J^\star_{\BB}$
when $\mu_0,\mu_1$ are compactly supported.

\begin{lemma}\label{lem:compact}
	Consider $c$ a real constant, $K\subset M$ a compact subset  and $\mu_0,\mu_1\in \mc P_2(M)$
	with compact supports included in $K$,
	then the set 
	\[
	\mc N(c,K):=\{\eta:\ J^\star_{\BB}(\eta)\leq c,\ \pi_\#\eta(I\times K)=1\}
	\]
	is compact in $\Y^{\mu_1}_{\mu_0}(I;HM)$.
\end{lemma}
\proof

We consider a point $x_0\in M$ and the distance function $x_0\mapsto d(x,x_0)$.
Then, we consider the integral 
\begin{equation}\label{boundtight}
\int_{I\times HM}d(x_0,\pi(v))+\norm{v}^2\,\d\eta(t,v).
\end{equation}
We want to show that it has bounded value
for $\eta\in \mc N(c,K)$.
The first term is constant on any fiber of the map $\pi\colon HM\to M$, therefore
\begin{equation}
\int_{I\times HM}d(x_0,\pi(v))\,\d\eta(t,v)=\int_{I\times M}d(x_0,x)\,\d\mu_t\d t,
\end{equation}
which is finite because $\pi_\#\eta=\mu$ is concentrated on a compact set
by hypothesis. 
Furthermore, the term $\int \norm{v}^2\,\d\eta$ equals $J^\star_{\BB}(\eta)$ and is also bounded by hypothesis.

Therefore for every $\eta\in \mc N(c,K)$ the integral \eqref{boundtight}
is bounded by a constant depending on~$x_0$, $\supp(\mu_0)$, $\supp(\mu_1)$ and $c$.
Moreover
 the integrand function 
\[v\mapsto d(\pi(v),x_0)+\norm{v}^2\]
has compact
level sets in $HM$.
By \cite[Corollary 3.61]{flogo12}, in the case of sets of Young measures
this implies the thightness of $\mc N(c,K)$ (see also \cite[Remark 3.53-(i)]{flogo12}), which in turn implies 
the
compactness of $\mc N(c,K)$ because $\mc P_2(M)$ is Polish with respect to the Wasserstein distance
for any Polish space $(M,d)$, and thigtness is equivalent to compactness on Polish spaces.\fine

\bibliography{mybiblio2}{}
\bibliographystyle{siam}

\end{document}